\documentclass[12pt,reqno,a4paper]{amsart}
\usepackage[utf8]{inputenc}

\usepackage{amsmath,amssymb,amsfonts,amsthm,amscd,mathtools,mathrsfs}
\usepackage{latexsym}
\usepackage{dsfont}
\usepackage{xcolor}
\usepackage{ytableau}
\usepackage[enableskew]{youngtab}
\usepackage{accents,upgreek}
\usepackage{pifont}

\usepackage{enumitem}
\usepackage{multicol,multirow,makecell}
\usepackage{nicematrix}
\usepackage{caption,subcaption}
\usepackage{setspace}

\usepackage{tikz}
\usetikzlibrary{matrix,arrows,decorations.pathmorphing,positioning}
\usepackage{pgfplots}
\pgfplotsset{compat=1.18}

\usepackage{verbatim}
\usepackage{soul}
\usepackage{dirtytalk}

\usepackage[letterpaper, margin=1in]{geometry}

\usepackage[
    colorlinks=true,
    linkcolor=blue,
    citecolor=green,
    urlcolor=blue
]{hyperref}
\usepackage[nameinlink,noabbrev,capitalize]{cleveref}
 
\theoremstyle{plain}
\newtheorem{thm}{Theorem}[section]

\newtheorem{lem}[thm]{Lemma}       
\newtheorem{cor}[thm]{Corollary}

\theoremstyle{definition}

\theoremstyle{remark}
\newtheorem{rem}[thm]{Remark}

\numberwithin{equation}{section}

\newcommand{\br}{\mathbb{R}}
\newcommand{\bc}{\mathbb{C}}
\newcommand{\bz}{\mathbb{Z}}
\newcommand{\bq}{\mathbb{Q}}
\newcommand{\bn}{\mathbb{N}}
\newcommand{\mf}{\mathfrak}

\newcommand{\mm}{\mathrm}

\makeatletter
\def\section{\def\@secnumfont{\mdseries}\@startsection{section}{1}%
    \z@{.7\linespacing\@plus\linespacing}{.5\linespacing}%
    {\normalfont\scshape\centering}}

\def\subsection{\def\@secnumfont{\bfseries}\@startsection{subsection}{2}%
    {\parindent}{.5\linespacing\@plus.7\linespacing}{-.5em}%
    {\normalfont\bfseries}}
\makeatother

\title{Generalized Chromatic Polynomials of Cycles and uni-cyclic Graphs}
\author{Dipnit Biswas}
\address{Indian Institute of Technology Madras, Chennai, India}
\email{ma25r011@smail.iitm.ac.in}
\thanks{}

\author{R. Venkatesh}
\address{Department of Mathematics, Indian Institute of Science, Bangalore 560012, India}
\email{rvenkat@iisc.ac.in}

\begin{document}

\begin{abstract} In \cite{R81}, Read obtained an explicit formula for the generalized chromatic polynomials of cycle graphs using the classical deletion–contraction principle. In this paper, we take a linear-algebraic approach and present a new proof of Read's formula by realizing the generalized chromatic polynomials of cycle graphs as traces of certain linear operators. As an application of this approach, we also derive explicit formulas for the generalized chromatic polynomials of uni-cyclic graphs.
\end{abstract}

\maketitle

\section{Introduction}

\noindent
Let $G$ be a simple graph (without any loops and multiple edges) with a vertex set $V=\{1,\dots,n\}$ and an edge set $E$. For $\mf{m} = (m_i)_{i=1}^n \in \bz_{\ge 0}^n$ and $q\in \bn$, a \emph{proper $q$-multi-coloring} of $G$ with weight vector $\mf{m}$ is a mapping $c:V\to P(\{1, \dots, q\})$
assigning to each $v\in V$ an $m_v$-element subset $c(v) \subseteq \{1,\dots,q\}$ such that $$\text{ $c(u)\cap c(v)=\emptyset$ whenever $uv\in E$.}$$
We say two such colorings $c_1$ and $c_2$ are distinct if there exists $v\in V$ such that $c_1(v)\neq c_2(v).$  For given $q\ge 1$, 
the total number of distinct such multi-colorings is a polynomial in $q$, denoted by $\pi^{\mf{m}}_G(q)$. This polynomial has a close relationship with the usual chromatic polynomials of graphs, see \cite[Remark 3.3]{ADV18}. Indeed, this generalization provides a uniform framework for studying the chromatic polynomials of the entire family of clan graphs (or clique blow-up graphs) associated to $G$.
The clan graph $G(\mf m)$ associated to $G$ and $\mf m$ is constructed as follows: for each $v$ with $m_v\neq 0$, take a clique (complete graph) of size $m_v$ with vertex set $\{v_1, \ldots, v_{m_v}\}$ and join all vertices of the $u$-th and $v$-th cliques if $uv \in E$. Then, we have
$$\pi^{\mf{m}}_G(q) = \frac{\pi_{G(\mf{m})}(q)}{\prod_{v=1}^{n} m_v!}.$$

\noindent
For instance, when \(G\) is chordal, these polynomials can be computed uniformly for all \(\mf m\); see \cite{BHV25}. However, the computation becomes considerably more complicated when \(G\) is not chordal. In this direction, the following theorem, due to Read \cite[Theorem 3]{R81}, gives an explicit formula for \(\pi^{\mf{m}}_G(q)\) in the first non-chordal case, namely, when \(G\) is a cycle of length more than $3$.
\begin{thm}\label{proper q-multi-coloring formula for C_n}
Let $n\in \bn$ with $n\ge 3$, and let $C_n$ be the cycle on the vertices $\{1,\dots,n\}$. Fix a weight vector $\mf{m} = (m_1,\dots,m_n)\in \bn^n$, set $m_{n+1} = m_1$ and define $m_* \coloneqq \min\{m_i \mid 1\le i \le n\}$. For every $q\in \bc$, we have
\begin{equation}\label{falling factorial formula for q-multi-coloring of C_n}
\pi^{\mf{m}}_{C_n}(q) = \frac{\prod_{i\,=\,1}^n (q)_{m_i + m_{i+1}}}{\prod_{i\,=\,1}^n m_i!} \sum_{k\,=\,0}^{m_*}(-1)^{kn}v_k(q) \prod_{i\,=\,1}^n \frac{(m_i)_k}{(q)_{m_i+k}}
\end{equation}
where $v_k(q) \coloneqq \binom{q}{k} - \binom{q}{k-1}$ for all $k\ge 0$, and $(x)_j \coloneqq x(x-1)(x-2)\cdots (x-j+1)$ denotes the falling factorial for any $x\in \br$ and $j\in \bn$ (with $(x)_0 \coloneqq 1$). \qed
\end{thm}

\noindent
The main goal of this article is to provide a different proof of this Theorem \ref{proper q-multi-coloring formula for C_n} by realizing the chromatic polynomials $\pi^{\mf{m}}_{C_n}(q)$ as traces of some natural linear operators.  In \cite{B01, B02, B04}, Biggs considered a similar operator and used the representation theory of the symmetric groups to compute the trace of this operator to give explicit formulas for bracelet graphs-which, in our setting, corresponds to clan graph of a cycle with diagonal weight vector $(m,\dots,m)$. In contrast, our proof bypasses the symmetric group representations and provides a considerably simpler and more direct derivation.
Using Read's result, we can also obtain explicit formulas for generalized chromatic polynomials of uni-cyclic graphs as an immediate corollary. 
\begin{thm}\label{proper q-multi-coloring formula for uni-cyclic graphs}
Let $G$ be a finite connected uni-cyclic graph whose unique cycle is $C_n=(c_1,c_2,\dots,c_n)$. For every $v\notin V(C_n)$, let $p(v)$ denote the unique neighbor of $v$ along the path from $v$ to $C_n$. Given a weight vector $\mf{m} = (m_v\mid v\in V)$, set
\[
a_i \coloneqq m_{c_i} \quad \text{for } 1\le i\le n,\qquad a_{n+1} \coloneqq a_1,\qquad a_* \coloneqq \min_{1\le i\le n} a_i.
\]
Then, for every $q\in \bc$,
\begin{equation}\label{formula for uni-cyclic graphs}
\pi_G^{\mf{m}}(q) = \pi_{C_n}^{(a_1,\dots,a_n)}(q) \prod_{v\notin V(C_n)} \binom{q - m_{p(v)}}{m_v}. \qed
\end{equation} 
\end{thm}

\section{Subset-function spaces and disjointness operator} \label{subset-function spaces}

\noindent
We recall some standard results from \cite[Chapter 4]{S13} in this section, and include some of their proofs for the sake of completeness.
For $q\ge 1$, we denote $[q]=\{1, \dots, q\}$. For $0\le r\le q$, let $\binom{[q]}{r}$ denote the set of all $r$-element subsets of $\{1,\dots,q\}$. Let $V_r$ be the space of $\bq$-valued functions on 
$\binom{[q]}{r}$:
\[
V_r \coloneqq \left\{ f \ \Big| \ f: \binom{[q]}{r} \to \bq \right\}.
\] It is a vector space with respect to pointwise additions and scalar multiplications of dimension $\binom{q}{r}$ over $\bq$.
For each $A\in \binom{[q]}{r}$, define the indicator function $1_A:\binom{[q]}{r} \to \bq$ by
\[
1_A(B) = 
\begin{cases}
    1 & \text{if } B = A,\\
    0 & \text{otherwise}.
\end{cases}
\]
The set $\mathcal{B}_r\coloneqq\{1_A\mid A\in \binom{[q]}{r}\} \subseteq V_r$ forms the standard basis for $V_r$. Define 
\[\langle f,g\rangle_{V_r} = \sum_{A\,\in\,\binom{[q]}{r}} f(A)g(A), \ \text{for all $f,g\in V_r$. }\] 
\subsection{}
For $r\ge 1$, the \emph{down-operator} $D_r : V_r \to V_{r-1}$ is a linear map defined by 
\[
(D_r(f))(A) \coloneqq \sum_{\substack{B\,\supseteq\,A \\|B|\,=\,r}} f(B),
\]
for all $A\in \binom{[q]}{r-1}$ and $f\in V_{r}$. Set $H_r \coloneqq \ker D_r$, and put $H_0 \coloneqq V_0$. 
For $r\ge 0$, the \emph{up-operator} $U_r: V_r \to V_{r+1}$ is a linear map defined by
\[
(U_r(f))(A) \coloneqq \sum_{\substack{B\,\subseteq\,A \\|B|\,=\,r}} f(B),
\]
for all $A\in \binom{[q]}{r+1}$ and $f\in V_{r}$.  The operator $U_{r-1}$ is the formal adjoint of $D_r$ with respect to this inner product (see \cite[Equation 4.4]{S13}). 
\begin{lem}\label{U_(r-1) is adjoint of D_r}
If $1\le r\le q$, then
\[
\langle D_r(f),g\rangle_{V_{r-1}} = \langle f, U_{r-1}(g)\rangle_{V_r}
\]
for all $f\in V_r$ and $g\in V_{r-1}$. \qed
\end{lem}
\noindent
We have the following commutation relation between down and up operators (see \cite[Lemma 4.6]{S13}). 
\begin{lem}\label{commutation}
For $0 \le j \le q-1$,
\[
D_{j+1} \circ U_j - U_{j-1} \circ D_j = (q-2j)\,\mathrm{id}_{V_j},
\]
where $U_{-1} \circ D_0 := 0$. \qed
\end{lem} \noindent It is easy to compute the dimension of $H_r$. 
The following Lemma also follows from \cite[Theorem 4.7]{S13}.
\begin{lem}\label{dimension of H_r}
If $0 \le r \le (q+1)/2$, then
\[
\dim H_r = \binom{q}{r} - \binom{q}{r-1},
\]
where $\binom{q}{-1} = 0$.
\end{lem}
\begin{proof}
If $r = 0$ there is nothing to prove, since $H_0 = V_0$ and $\dim V_0 = 1$.
Let $1 \le r \le (q+1)/2$ and $g \in V_{r-1}$. By Lemma~\ref{commutation} with
$j = r-1$, and Lemma~\ref{U_(r-1) is adjoint of D_r},
\begin{align*}
\langle D_r U_{r-1}(g), g \rangle_{V_{r-1}}
&= \langle U_{r-2} D_{r-1}(g), g \rangle_{V_{r-1}} + (q-2r+2)\,\langle g, g \rangle_{V_{r-1}} \\
&= \langle D_{r-1}(g), D_{r-1}(g) \rangle_{V_{r-2}} + (q-2r+2)\,\langle g, g \rangle_{V_{r-1}},
\end{align*}
where for $r = 1$ the first term on the right is absent. Since $2r \le q+1$, we have
$q-2r+2 \ge 1$, so the right-hand side is positive whenever $g \ne 0$. Hence
$D_r \circ U_{r-1}$ is injective on $V_{r-1}$, and therefore bijective since
$V_{r-1}$ is finite-dimensional. In particular $D_r$ is surjective, so
$\operatorname{rank}(D_r) = \dim V_{r-1} = \binom{q}{r-1}$. By the rank--nullity theorem,
\[
\dim H_r = \dim \ker D_r = \dim V_r - \operatorname{rank}(D_r)
= \binom{q}{r} - \binom{q}{r-1}. \qedhere
\]
\end{proof}

\noindent
The next Lemma is immediate from the above argument (see \cite[Theorem 4.7]{S13}). 
\begin{lem}\label{injectivity of U_r}
For any $0\le r < q/2$, the up-operator $U_r:V_r\to V_{r+1}$ is injective. \qed
\end{lem}

\subsection{}
It is easy to see that the space $V_r$ admits a decomposition in terms of up-operators and kernel of down-operators.
\begin{lem}\label{decomposition of V_r}
If $0\le r\le q/2$, then
\[
V_r = H_r \oplus \bigoplus_{k\,=\,1}^r U_{r-1}\circ U_{r-2} \circ \cdots \circ U_{r-k} (H_{r-k}).
\]
\end{lem}

\begin{proof}
The case $r=0$ is vacuously true. Let $r\ge 1$. By Lemma \ref{U_(r-1) is adjoint of D_r}, $U_{r-1}$ is the adjoint of $D_r$. Thus by duality of standard inner product, we have 
$U_{r-1}(V_{r-1})^{\perp}=H_r$. Since the standard inner product on $V_r$ is
positive definite, we get
\[
V_r = H_r \oplus U_{r-1}(V_{r-1}).
\]
 Now the result follows by induction on $r$.
\end{proof}

\subsection{}
The next result describes the composition of successive down-operators.
\begin{lem}\label{composition of successive down-operators}
Let $0\le k < r \le q$. Then for any $f\in V_r$,
\[
(D_{k+1} \circ D_{k+2}\circ \cdots \circ D_r (f))(A) = (r-k)!\sum_{\substack{B \,\supseteq\, A\\ |B|\, =\, r}} f(B)
\]
for all $A\in \binom{[q]}{k}$.
\end{lem}

\begin{proof}
Fix $f\in V_r$ and $A\in \binom{[q]}{k}$. Then 
\begin{align*}
(D_{k+1} \circ D_{k+2}\circ \cdots \circ D_r (f))(A) & = \sum_{\substack{A_1\,\supseteq\, A\\|A_1|\,=\,k+1}} (D_{k+2}\circ \cdots \circ D_r (f))(A_1) \\
& = \sum_{\substack{A_1\,\supseteq\, A\\|A_1|\,=\,k+1}} \sum_{\substack{A_2\,\supseteq\, A_1\\|A_2|\,=\,k+2}} (D_{k+3}\circ \cdots \circ D_r (f))(A_2) \\
& = \sum_{\substack{A_1\,\supseteq\, A\\|A_1|\,=\,k+1}} \sum_{\substack{A_2\,\supseteq\, A_1\\|A_2|\,=\,k+2}} \cdots \sum_{\substack{A_{r-k}\,\supseteq\, A_{r-k-1}\\|A_{r-k}|\,=\,r}} f(A_{r-k}) \\
& = C_{r,k} \sum_{\substack{A_{r-k}\,\supseteq\, A\\|A_{r-k}|\,=\,r}} f(A_{r-k})
\end{align*}
where $C_{r,k}$ is the number of chain of subsets
\[
A \subseteq A_1 \subseteq A_2 \subseteq \cdots \subseteq A_{r-k-1} \subseteq A_{r-k}.
\]
such that $|A_i| = k+i$ for all $1\le i \le r-k-1$. Such a chain is determined by the order
in which the $r-k$ elements of $A_{r-k}\setminus A$ are added, so $C_{r,k} = (r-k)!$ for
every $A_{r-k}$. This completes the proof.
\end{proof} \noindent
We have the following immediate corollary. 
\begin{cor}\label{action of f in H_r on k<r element subsets}
Fix $0\le k< r\le q$. If $f\in H_r$, then
\[
\sum_{\substack{B \,\supseteq\, A\\ |B|\,=\,r}} f(B) = 0
\]
for all $A\in \binom{[q]}{k}$. \qed
\end{cor}
\noindent 
By an argument similar to that of Lemma~\ref{composition of successive down-operators}
(or by taking adjoints in it, using Lemma~\ref{U_(r-1) is adjoint of D_r}), the
composition of successive up-operators is described as follows: for
$0\le k< r\le q$ and any $f\in V_k$,
\[
(U_{r-1}\circ U_{r-2}\circ \cdots \circ U_k(f))(A) = (r-k)!\sum_{\substack{B \,\subseteq\, A\\ |B|\, =\, k}} f(B)
\qquad \text{for all } A\in \tbinom{[q]}{r}.
\]
Motivated by the above expression, for $0\le k\le r\le q$, we define
$L_{k,r}: V_k \to V_r$ by
\[
(L_{k,r}(f))(A) = \sum_{\substack{B \,\subseteq\, A\\ |B|\, =\, k}} f(B)
\]
for all $f\in V_k$ and $A\in \binom{[q]}{r}$. Thus $L_{r,r}$ is the identity operator,
and $$\text{ $U_{r-1}\circ U_{r-2}\circ \cdots \circ U_k = (r-k)!\,L_{k,r}$ for $k<r$.}$$

\subsection{}
For non-negative integers $a,b$ such that $a+b\le q$, define the \emph{disjointness operator} $T_{a,b}: V_b \to V_a$ by
\[
(T_{a,b}(f))(A) = \sum_{\substack{|B| \,=\, b \\ A\cap B\,=\, \emptyset}} f(B) \ \ \text{for all $A\in \binom{[q]}{a}$.}
\]
\noindent
We need the following lemma.
\begin{lem}\label{action of disjointness operator on up-operator}
Let $a,b\ge 0$ with $a+b\le q$, and let $0\le r \le \min\{a,b\}$. For any $f\in H_r$,
\[
(T_{a,b}\circ L_{r,b})(f) = (-1)^r \binom{q-a-r}{b-r}L_{r,a}(f).
\]
\end{lem}

\begin{proof}
Fix $f\in H_r$ and $A\in \binom{[q]}{a}$. Then 
\begin{align*}
(T_{a,b}\circ L_{r,b}(f))(A) & = \sum_{\substack{|B|\,=\,b\\A\cap B\,=\,\emptyset}} (L_{r,b}(f))(B)  = \sum_{\substack{|B|\,=\,b\\A\cap B\,=\,\emptyset}} \sum_{\substack{C\,\subseteq \,B \\ |C|=r}} f(C) = \binom{q-a-r}{b-r}  \sum_{\substack{|C|\,=\,r\\ C\cap A\,=\,\emptyset}} f(C).
\end{align*}
since for each $r$-set $C$ disjoint from $A$, the $b$-sets $B$ with $C\subseteq B$ and
$A\cap B=\emptyset$ are obtained by choosing $b-r$ elements of $[q]\setminus(A\cup C)$,
a set of size $q-a-r$. For $0\le i\le r$, set
\[
F_i = \sum_{\substack{|C|\,=\,r\\ |C\cap A|\,=\,i}} f(C).
\]
Then by applying Corollary \ref{action of f in H_r on k<r element subsets}, we get
\[
\sum_{i\,=\,t}^r \binom{i}{t} F_i = \sum_{i\,=\,t}^r \sum_{\substack{|C|\,=\,r\\ |C\cap A|\,=\, i}} \binom{i}{t} f(C) = \sum_{\substack{B\,\subseteq\,A\\|B| \,=\, t}} \sum_{\substack{C\,\supseteq\,B\\|C| \,=\, r}} f(C) = 0 
\]
for all $0\le t < r$. This implies
\[
\begin{pmatrix}
    \binom{0}{0} & \binom{1}{0} & \cdots & \binom{r-1}{0} & \binom{r}{0} \vspace{1mm} \\ 
    0 & \binom{1}{1} & \cdots & \binom{r-1}{1} & \binom{r}{1} \\
    \vdots & \vdots & \ddots & \vdots & \vdots \\
    0 & 0 & \cdots & \binom{r-1}{r-1} & \binom{r}{r-1}
\end{pmatrix}
\begin{pmatrix}
    F_0 \\
    F_1 \\
    \vdots \\
    F_r
\end{pmatrix} =
\begin{pmatrix}
    0 \\
    0 \\
    \vdots \\
    0
\end{pmatrix}.
\]
Solving this system of equations, we get
\[
F_t = (-1)^{r-t}\binom{r}{t} F_r
\]
for all $0\le t < r$. Thus $F_0 = (-1)^rF_r$. Hence
\[
\sum_{\substack{|C|\,=\,r\\ C\cap A\,=\,\emptyset}} f(C) = \sum_{\substack{|C|\,=\,r\\ |C\cap A|\,=\,0}} f(C) = (-1)^r \sum_{\substack{|C|\,=\,r\\ |C\cap A|\,=\,r}} f(C) = (-1)^r\sum_{\substack{C\,\subseteq\,A\\|C|\,=\,r}} f(C) = (-1)^r(L_{r,a}(f))(A).
\]
This completes the proof.
\end{proof}

\section{Main Theorem} \noindent
We are ready to provide a new proof of Read's result. 
\begin{thm}
Let $n\in \bn$ with $n\ge 3$, and let $C_n$ be the cycle graph on the vertex $\{1,\dots,n\}$. Fix a weight vector $\mf{m} = (m_1,\dots,m_n)\in \bn^n$, set $m_{n+1} = m_1$ and define $m_* \coloneqq \min\{m_i \mid 1\le i \le n\}$. As polynomials in $q$, we have
\begin{equation}\label{binomial formula for q-multi-coloring of C_n}
\pi^{\mf{m}}_{C_n}(q) = \sum_{k\,=\,0}^{m_*}(-1)^{kn} v_k(q)\prod_{i\,=\,1}^n \binom{q-m_i-k}{m_{i+1}-k}
\end{equation}
where $v_k(q) = \binom{q}{k} - \binom{q}{k-1}$.
\end{thm}

\begin{proof} Both sides of the formula are polynomials in $q$, so it suffices to prove it for
every integer $q \ge 2\max\{m_i \mid 1\le i \le n\}$; fix such a $q$. Both sides of Equation \ref{binomial formula for q-multi-coloring of C_n} are also unchanged
when the vertices of $C_n$ are relabelled cyclically, so we may assume that
$m_1 = m^*$. 
Consider the operator
\[
P = T_{m_1,m_2}\circ T_{m_2,m_3}\circ \cdots \circ T_{m_n,m_1} : V_{m_1} \to V_{m_1},
\]
it is well-defined since $m_i + m_{i+1} \le q$. We compute the trace of this operator, $P$, in two different ways. Computing its trace with respect to the standard basis gives us the generalized chromatic polynomials of $C_n$, that is
 $\mathrm{Trace}(P) = \pi^{\mathbf m}_{C_n}(q)$.
For $B \in \binom{[q]}{b}$ we have
$$T_{a,b}(1_B) = \sum_{|A|\, =\, a,\ A\cap B\, = \,\emptyset} 1_A,$$ so the $(A,B)$ entry of the
matrix of $T_{a,b}$ with respect to $\mathcal{B}_b$ and $\mathcal{B}_a$ is $1$ if
$A \cap B = \emptyset$ and $0$ otherwise. Multiplying these matrices, the diagonal
entry of the matrix of $P$ indexed by $A_1 \in \binom{[q]}{m_1}$ is the number of
sequences $$\text{$(A_2,\dots,A_n)$ with $A_i \in \binom{[q]}{m_i}$ and
$A_1\cap A_2 = A_2 \cap A_3 = \cdots = A_n \cap A_1 = \emptyset$.}$$ Summing over $A_1$,
$\mathrm{Trace}(P)$ is the number of sequences $(A_1,\dots,A_n)$ with
$A_i \in \binom{[q]}{m_i}$ and $A_i \cap A_{i+1} = \emptyset$ for $1 \le i \le n$
(where $A_{n+1} = A_1$), which is $\pi^{\mathbf m}_{C_n}(q)$.

\medskip
\noindent On the other hand, computing the $\mathrm{Trace}(P)$ using Lemma \ref{action of disjointness operator on up-operator} gives us the right hand side of the Equation \ref{binomial formula for q-multi-coloring of C_n}.
By Lemma \ref{decomposition of V_r} and definition of $L_{k,r}$, we have
\begin{align*}
V_{m_1} & = H_{m_1} \oplus \bigoplus_{k\,=\,1}^{m_1} U_{m_1-1}\circ U_{m_1-2} \circ \cdots \circ U_{m_1-k} (H_{m_1-k}) 
= \bigoplus_{k\,=\,0}^{m_1} L_{k,m_1}(H_k).
\end{align*}
Fix $0\le k\le m_1$. Then by Lemma \ref{action of disjointness operator on up-operator}, we have
\begin{align*}
& \hspace{1mm} T_{m_1,m_2}\circ T_{m_2,m_3}\circ \cdots\circ T_{m_{n-1},m_n}\circ T_{m_n,m_1}\circ L_{k,m_1}(f)\\ = & \hspace{1mm} (-1)^k\binom{q-m_n - k}{m_1 - k} T_{m_1,m_2}\circ T_{m_2,m_3}\circ \cdots\circ T_{m_{n-1},m_n}\circ L_{k,m_n}(f)\\ =
& \hspace{1mm} (-1)^{2k}\binom{q-m_n - k}{m_1 - k}\binom{q-m_{n-1}-k}{m_n-k} T_{m_1,m_2}\circ T_{m_2,m_3}\circ \cdots\circ T_{m_{n-2},m_{n-1}}\circ L_{k,m_{n-1}}(f)\\ =
& \hspace{1mm} (-1)^{kn}\prod_{i\,=\,1}^n \binom{q - m_i - k}{m_{i+1}-k} L_{k,m_1}(f)
\end{align*}
for all $f\in H_k$. Since $U_r$ is injective for all $r\ge 0$ (by Lemma \ref{injectivity of U_r}), $L_{k,m_1}$ is injective for all $0\le k\le m_1$. Thus $\dim L_{k,m_1}(H_k) = \dim H_k = \binom{q}{k} - \binom{q}{k-1}$ (by Lemma \ref{dimension of H_r}). Hence,
\begin{align*}
 \mathrm{Trace}(T_{m_1,m_2}\circ T_{m_2,m_3} \circ \cdots\circ T_{m_n,m_1}) 
& = \sum_{k=0}^{m_1} (-1)^{kn}\prod_{i\,=\,1}^n \binom{q - m_i - k}{m_{i+1}-k} \left[\binom{q}{k} - \binom{q}{k-1}\right] \\ 
& = \sum_{k=0}^{m_*} (-1)^{kn} v_k(q) \prod_{i\,=\,1}^n \binom{q - m_i - k}{m_{i+1}-k}.
\end{align*}
This completes the proof. 
\end{proof}
\begin{rem}
For indeterminate $q$, substituting 
\[
\binom{q - m_i -k}{m_{i+1} - k} = \frac{(q)_{m_i + m_{i+1}}}{(q)_{m_i+k}}\frac{(m_{i+1})_k}{m_{i+1}!}
\]
in the identity given by Equation \ref{binomial formula for q-multi-coloring of C_n} we get the identity given by Equation \ref{falling factorial formula for q-multi-coloring of C_n}. Hence Theorem \ref{proper q-multi-coloring formula for C_n} now follows.
\end{rem}

\subsection{Uni-cyclic graphs} \label{uni-cyclic graphs}
We now provide proof of Theorem \ref{proper q-multi-coloring formula for uni-cyclic graphs}.
First, we assume that $q$ is sufficiently large. The cycle $C_n = (a_1,\dots,a_n)$ has $\pi_{C_n}^{(a_1,\dots,a_n)}(q)$ many proper $q$-multi-colorings. Fix one such proper $q$-multi-coloring. Next we color the remaining vertices of $G$ starting with the vertices which are at a distance of $1$ from the cycle and proceed iteratively in increasing order of distance. 

Therefore when we color a vertex $v$ of $G$ outside the cycle $C_n$, the only neighbor of $v$ which is already colored is $p(v)$. Thus, the number of possible colorings of $v$ is 
\[
\binom{q-m_{p(v)}}{m_v}.
\]
Hence, taking the product of $\pi_{C_n}^{(a_1,\dots,a_n)}(q)$ with the possible proper $q$-multi-colorings of vertices of $G$ not on $C_n$, gives us Equation \ref{formula for uni-cyclic graphs}. Since both sides of Equation \ref{formula for uni-cyclic graphs} are polynomials in $q$, and they give the same value when evaluated at any sufficiently large $q$, so the identity holds for all $q\in \bc$.

\section*{Declaration of AI use}

\noindent The authors used Google Gemini (Gemini 3.6 Flash) and Claude (Opus 5.5) to improve readability of the prose, refine technical phrasing, and correct for grammatical errors. The authors thoroughly reviewed and revised all results and take full responsibility for the final content of the article.

\bibliographystyle{plain}
\bibliography{bibliography}

\end{document}